\documentclass[12pt]{article}
\usepackage{graphics}
\usepackage{graphicx}
\usepackage{color}
\newtheorem{remark}{Remark}[section]

\title{A mathematical model for the atomic clock error in case of jumps}

\author{C Zucca
\\
\small Department of Mathematics \lq\lq G. Peano\rq\rq,\\ 
\small U\-ni\-ver\-si\-ty of To\-ri\-no,\\ 
\small Via Carlo Al\-ber\-to 10, 10123 Turin, Italy\\
\small \texttt{cristina.zucca@unito.it}\\
\and
P Tavella\\
\small Istituto Nazionale di Ricerca Metrologica (INRiM),\\ 
\small Strada delle Cacce 91,\\ 
\small 10135 Turin, Italy}

\date{}
\begin{document}
\maketitle

\begin{abstract}
We extend the mathematical model based on stochastic differential equations describing the error gained by an atomic clock to the cases of anomalous behavior including jumps and an increase of instability. We prove an exact iterative solution that can be useful for clock simulation, prediction, and interpretation, as well as for the understanding of the impact of clock error in the overall system in which clocks may be inserted as, for example, the Global Satellite Navigation Systems.
\end{abstract}

\vspace{2pc}

\section{Introduction}
Atomic clocks are usually the heart of complex scientific and technological systems and clock behavior has a direct impact on the overall performance of the system. For the study, simulation, analysis, characterization, and interpretation of all these systems it is therefore mandatory to have a consistent mathematical clock model able to represent its typical behavior. Several studies have dealt with this matter \cite{A, AD, G, GSTZ,VB, VDBT, ZT}.

Recently, particularly in space applications such as the Global Satellite Navigation Systems (GNSS), anomalous clock behaviors have been experimentally observed \cite{BHS,C,CFSGT,DBWREW,GT,GT1,W}: phase or frequency jumps, increase of instability, general non stationary trends. Such anomalous behaviors have an important impact on the global system and may affect the overall performance very seriously. This could have, in case of GNSS, dramatic consequences even security of life.

In this paper we propose a mathematical model including the typical anomalies of atomic clocks and we show how the stochastic multidimensional differential equation describing such behaviors can be exactly solved obtaining the complete description of clock states. In addition, the solution can be written in an iterative form. 

We therefore provide a clock model including anomalies with an exact iterative form that may be useful for simulations as well as for the insertion of the clock model in further processing as Kalman filtering, orbit and clock estimation algorithms, time scale algorithms, or other types of estimation and filtering. Moreover we provide an exact solution that allows the prediction and evaluation of the clock error at a certain time after synchronization also in presence of anomalies and as a function of the epoch of the anomaly occurrence. This knowledge is of fundamental importance for GNSS evaluations as the  clock time deviation directly impact the user solution accuracy \cite{PT,T}.

In Section 2 we recall the main features of the mathematical model based on stochastic differential equations including phase, frequency, and drift jumps and we introduce in the model a possible sequence of jumps and the increase of instability.
In Section 3 we illustrate through examples the features captured by the mathematical model. In Section 4 we explicitly address the example of a space Rubidium clock on board of a navigation satellite and we evaluate the prediction error at a certain time after synchronization in presence of a frequency jump.

To ease the use of this model in simulating clock behavior, we also add the availability of a Matlab code that can be freely downloaded from the INRIM web site.

\section{The model of the atomic clock error with anomalies}

It is experimentally observed that the signal of an atomic clock can be affected by several anomalies. Our aim is to extend the model proposed in \cite{ZT} introducing anomalies. In the forthcoming we will consider three different types of anomalies that involve jumps at deterministic times or changes in the variance.

\subsection{Instantaneous jumps}

Let us consider the three-state clock model that we proposed in \cite{ZT}. The three-state clock model ${\bf X}=\{(X_1,X_2,X_3)(t), t\geq 0\}$ is described by the three dimensional stochastic differential equation:
\begin{eqnarray}\label{X}
\left\{
\begin{array}{lll}
dX_1(t)=(X_2(t)+\mu_1)dt+\sigma_1 dW_1(t)\\
dX_2(t)=(X_3(t)+\mu_2)dt+\sigma_2 dW_2(t)\\
dX_3(t)=\mu_3 dt+\sigma_3 dW_3(t)
\end{array}
\right.
\end{eqnarray}
with initial condition ${\bf X}(0)=(c_1,c_2,c_3)$, where $\{W_i(t),t\geq 0\}$, $i=1,2,3$, are three indipendent, one-dimensional standard Wiener processes. 

The variable $X_1(t)$ represents the time deviation, the derivative $\dot{X}_1(t)$ represents the frequency deviation of which $X_2(t)$ is only a component (i.e., what is generally called the random walk component). Finally $X_3(t)$ represents the frequency drift or aging. 

The metrological interpretation of the model, the solution of (\ref{X}) and its iterative form can be found in \cite{ZT}. Here we recall that the probability distribution function of the solution $\bf X$ at time $t$ is Normal with mean \cite{ZT} 
\begin{eqnarray}\label{M}
{\bf M}(t)=\left[
\begin{array}{l}
c_1+(c_2 + \mu_1) t + (c_3 + \mu_2) \frac{t^2}{2}+ \mu_3 \frac{t^3}{6}\\
c_2 + (c_3 + \mu_2) t + \mu_3 \frac{t^2}{2}\\
c_3+\mu_3 t
\end{array}
\right]
\end{eqnarray}
and covariance matrix
\begin{equation}\label{Sigma}
{\bf \Sigma}(t)=\left[
\begin{array}{lll}
\sigma_1^2 t+\sigma_2^2 \frac{t^3}{3}+\sigma_3^2 \frac{t^5}{20} & \sigma_2^2 \frac{t^2}{2}+\sigma_3^2 \frac{t^4}{8} & \sigma_3^2 \frac{t^3}{6}\\
 \sigma_2^2 \frac{t^2}{2}+\sigma_3^2 \frac{t^4}{8} & \sigma_2^2 t+\sigma_3^2 \frac{t^3}{3}&\sigma_3^2 \frac{t^2}{2}\\
\sigma_3^2 \frac{t^3}{6} &\sigma_3^2 \frac{t^2}{2}& \sigma_3^2 t
\end{array}
\right].
\end{equation}

Let us extend the model of the clock error (\ref{X})
characterized by jumps with amplitude $a_i$, $i=1,2,3$ in each component: phase, frequency and drift. To introduce a simple and  treatable mathematical description, we make use of general continuous time Markov Chain  ${\bf N}(t)=\{(N_1,N_2,N_3)(t),t\geq 0\}$, i.e. a continuous-time stochastic process which takes values in some finite or countable set and for which the time spent in each state has an exponential distribution (cf. \cite{S}). If we introduce these anomalies into the clock model (\ref{X}) the corresponding mathematical model of the clock becomes
\begin{eqnarray}\label{X anomalies}
\left\{
\begin{array}{lll}
dX_1(t)=(X_2(t)+\mu_1)dt+\sigma_1 dW_1(t)+a_1 dN_1(t)\\
dX_2(t)=(X_3(t)+\mu_2)dt+\sigma_2 dW_2(t)+a_2 dN_3(t)\\
dX_3(t)=\mu_3 dt+\sigma_3 dW_3(t)+a_3 dN_3(t)
\end{array}
\right.,
\end{eqnarray}
where $N_i(t)$, $i=1,2,3$ are independent continuous time Markov Chains.

\begin{figure}[t!] 
\centering 
\includegraphics[scale=0.6]{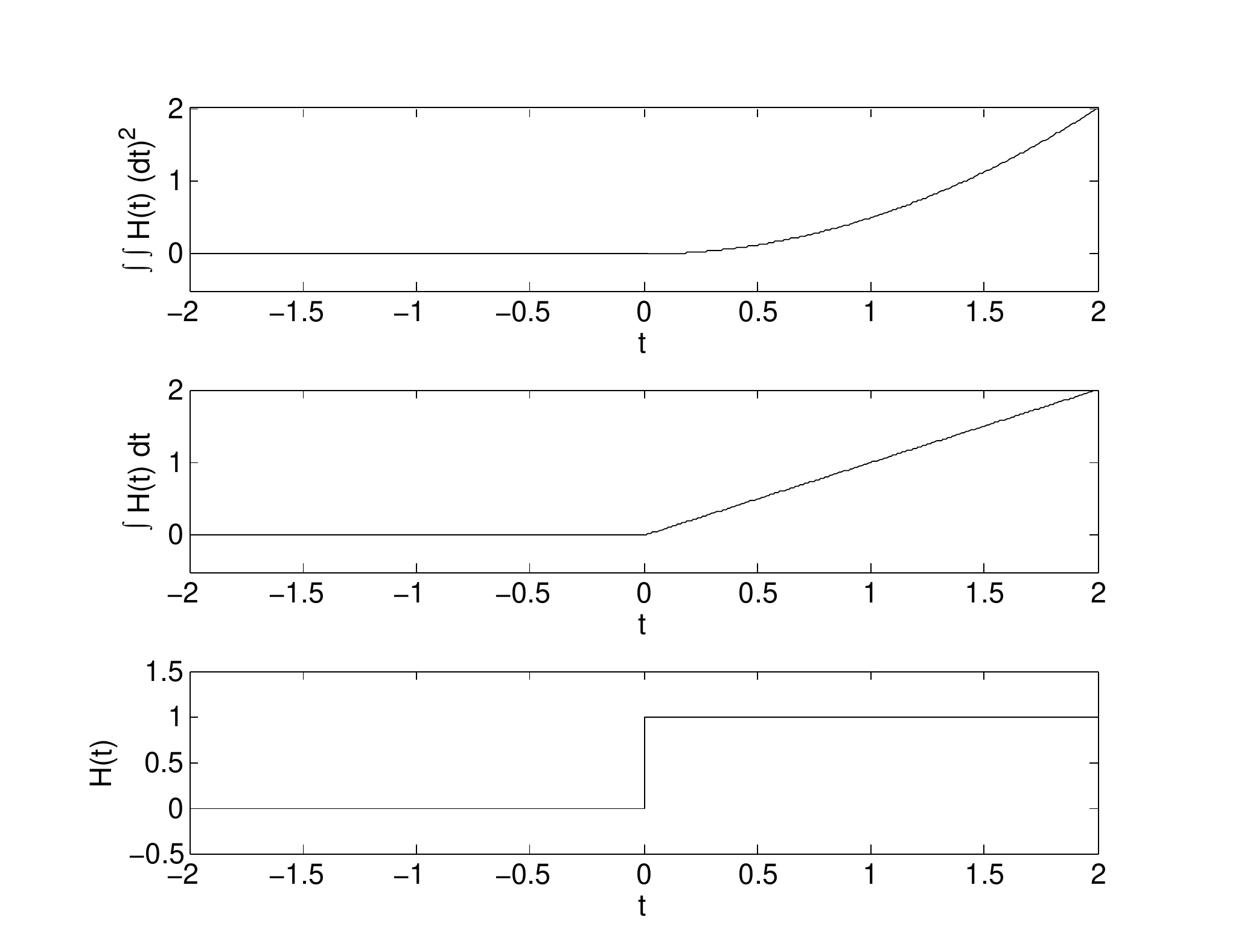}
\caption{Plot of the Heaviside function $H(t)$, its integral and its double integral (from below to above).}
\label{Fig:integrale H}
\end{figure}

The easiest particular case considers istantaneous jumps with amplitude $a_i$, $i=1,2,3$ that take place at deterministic times $\theta_i$, $i=1,2,3$. It means that 
\begin{equation}\label{N=H}
dN_i(t)=dH(t-\theta_i)=\delta(t-\theta_i)dt.
\end{equation}
where $H(t)$ is the Heaviside function and $\delta(t)$ is the Dirac delta function \cite{AS}.
Note that the derivatives of the Heaviside function is the Dirac delta function, the integral is the so called ramp function and the multiple integrals are given by
\begin{equation}\label{Int H}
\int \dots \int H(t) \underbrace{dt \dots dt}_n=\frac{t^n}{n!}H(t).
\end{equation}
In Figure \ref{Fig:integrale H} the Heaviside function, its integral and double integral are shown.

A first generalization of this simple model would consider $N_i(t)$, $i=1,2,3$ independent Poisson processes, i.e. anomalies occurring on random epochs, the intertimes between the anomalies are exponentially distributed and the jumps are deterministic of size 1 \cite{S}. 
This model will be the subject of a future study as it would help in treating the complete clock behaviour for example in the analysis of time scale algorithms. 

Here we consider the choice (\ref{N=H}), the system (\ref{X anomalies}) becomes
\begin{eqnarray}\label{X jump}
\left\{
\begin{array}{lll}
dX_1(t)=(X_2(t)+\mu_1)dt+\sigma_1 dW_1(t)+a_1 dH(t-\theta_1)\\
dX_2(t)=(X_3(t)+\mu_2)dt+\sigma_2 dW_2(t)+a_2 dH(t-\theta_2)\\
dX_3(t)=\mu_3 dt+\sigma_3 dW_3(t)+a_3 dH(t-\theta_3)
\end{array}
\right..
\end{eqnarray}
Integrating (\ref{X jump}) and recalling 
\begin{eqnarray}\label{delta}
\int_0^t a \delta(s-t_0)ds= aH(t-t_0) \nonumber \\
\int_0^t a H(s-t_0)ds= a(t-t_0)H(t-t_0)\\
\int_0^t a (s-t_0)H(s-t_0)ds= \frac{a(t-t_0)^2}{2}H(t-t_0),\nonumber
\end{eqnarray}
the solution can be written in closed form
\begin{eqnarray}\label{Xsol jump}
\left\{
\begin{array}{lll}
X_1(t)=c_1+(c_2+\mu_1)t+(c_3+\mu_2)\frac{t^2}{2}+\mu_3\frac{t^3}{6}+\sigma_1 W_1(t)\\
\quad \quad \quad \quad+\sigma_2\int_0^t(t-s)dW_2(s)+\sigma_3\int_0^t\frac{(t-s)^2}{2}dW_3(s)\\
\quad \quad \quad \quad +a_1H(t-\theta_1)+a_2(t-\theta_2)H(t-\theta_2)\\
\quad \quad \quad \quad+a_3\frac{(t-\theta_3)^2}{2}H(t-t_3)\\
X_2(t)=c_2+(c_3+\mu_2)t+\mu_3\frac{t^2}{2}+\sigma_2 W_2(t)+\sigma_3\int_0^t(t-s)dW_3(s)\\
\quad \quad \quad \quad +a_2H(t-\theta_2)+a_3(t-\theta_3)H(t-\theta_3)\\
X_3(t)=c_3+\mu_3 t+\sigma_3 W_3(t)+a_3H(t-\theta_3)
\end{array}
\right..
\end{eqnarray}

The probability distribution function of the solution $\bf X$ at time $t$ is Normal with mean ${\bf M^I}(t)$ and covariance matrix ${\bf \Sigma}(t)$,  where the mean depends on the mean value of the process without anomalies (\ref{M}) 
\begin{eqnarray}\label{M I}
\hspace{-1.5cm}{\bf M^I}(t)={\bf M}(t)+\left[
\begin{array}{l}
a_1H(t-\theta_1)+a_2(t-\theta_2)H(t-\theta_2)+a_3\frac{(t-\theta_3)^2}{2}H(t-t_3)\\
a_2H(t-\theta_2)+a_3(t-\theta_3)H(t-\theta_3)\\
a_3H(t-\theta_3)
\end{array}
\right]
\end{eqnarray}
while its covariance matrix (\ref{Sigma}) does not change.

Let us now consider a fixed time interval $[0,T]$ and an equally spaced partition $0\equiv t_0<t_1<\dots<t_N\equiv T$  and let us denote with $\tau=t_{k+1}-t_k$, $k=0,1,\dots,N-1$ the resulting discretization time step. For simplicity we hypothesize that the jumps can occur only in the discretization epochs. The sampling time can be as small as necessary to ensure this assumption. Moreover, observing that 
\begin{eqnarray}\label{H-H}
H(t_{k+1}-\theta)-H(t_k-\theta)=\delta_{(t_{k+1}-\theta)},
\end{eqnarray}
and consequently
\begin{eqnarray}\label{H-H-1}
(t_{k+1}-\theta)[H(t_{k+1}-\theta)-H(t_k-\theta)]=0,
\end{eqnarray}
where  $\delta_{(t_{k+1}-\theta)}$ is the Kronecker delta \cite{A},
we can express the solution at time $t_{k+1}$ in terms of the solution at time $t_k$ as:

\begin{eqnarray}\label{Xsol iter jump bis}
\left\{
\begin{array}{lll}
X_1(t_{k+1})=X_1(t_k)+(\mu_1+X_2(t_k))\tau+(\mu_2+X_3(t_k))\frac{\tau^2}{2}+\mu_3\frac{\tau^3}{6}+J_{k,1}\\
\quad \quad \quad \quad +a_1 \delta_{(t_{k+1}-\theta_1)}\\
X_2(t_{k+1})=X_2(t_k)+(\mu_2+X_3(t_k))\tau+\mu_3\frac{\tau^2}{2}+J_{k,2}+a_2 \delta_{(t_{k+1}-\theta_2)}\\
X_3(t_{k+1})=X_3(t_k)+\mu_3 \tau+J_{k,3}+a_3 \delta_{(t_{k+1}-\theta_3)}
\end{array}
\right.,
\end{eqnarray}
where ${\bf J}_k=(J_{k,1},J_{k,2},J_{k,3})$ can be interpreted as an innovation
\begin{equation}\label{innovation}
{\bf J}_k=\left[
\begin{array}{ll}
\sigma_1 (W_1(t_{k+1})-W_1(t_{k}))+\sigma_2\int_{t_k}^{t_{k+1}}(W_2(s)-W_2(t_k))ds\\
\quad +\sigma_3\int_{t_k}^{t_{k+1}}(t_{k+1}-s)(W_3(s)-W_3(t_k))ds\\
\sigma_2 (W_2(t_{k+1})-W_2(t_{k}))+\sigma_3\int_{t_k}^{t_{k+1}}(W_3(s)-W_3(t_k))ds\\
\sigma_3 (W_3(t_{k+1})-W_3(t_k))
\end{array}
\right].
\end{equation}
and is Normal distributed with zero mean and covariance matrix ${\bf Q}$ given by 
\begin{equation}\label{Q}
{\bf Q}=\left[
\begin{array}{lll}
\sigma_1^2 \tau+\sigma_2^2 \frac{\tau^3}{3}+\sigma_3^2 \frac{\tau^5}{20} & \sigma_2^2 \frac{\tau^2}{2}+\sigma_3^2 \frac{\tau^4}{8} & \sigma_3^2 \frac{\tau^3}{6}\\
 \sigma_2^2 \frac{\tau^2}{2}+\sigma_3^2 \frac{\tau^4}{8} & \sigma_2^2 \tau+\sigma_3^2 \frac{\tau^3}{3}&\sigma_3^2 \frac{\tau^2}{2}\\
\sigma_3^2 \frac{\tau^3}{6} &\sigma_3^2 \frac{\tau^2}{2}& \sigma_3^2h
\end{array}
\right].
\end{equation}

In the iterative form, we note that the jump acts only once on the corresponding component, then the iterative process carries the effect in the other components.

\subsection{Effect of two equal and opposite jumps}
Let us consider the model of the clock error 
characterized by an anomaly on the second component, a temporary frequency jump. The anomaly is composed of two istantaneous jumps at times $\theta_0$ and $\theta_1>\theta_0$. At time $\theta_0$ we register a positive jump of amplitude $a/\Delta$ and at time $\theta_1$ we register a negative jump of amplitude $-a/\Delta$, where $\Delta=\theta_1-\theta_0$.
The atomic clock error can be described by the following stochastic differential equation 
\begin{eqnarray}\label{X 2jump}
\hspace{-1cm}\left\{
\begin{array}{lll}
dX_1(t)=(X_2(t)+\mu_1)dt+\sigma_1 dW_1(t)\\
dX_2(t)=(X_3(t)+\mu_2)dt+\sigma_2 dW_2(t)+\frac{a}{\Delta} dH(t-\theta_0)-\frac{a}{\Delta} dH(t-\theta_1)\\
dX_3(t)=\mu_3 dt+\sigma_3 dW_3(t).
\end{array}
\right.
\end{eqnarray}
Integrating (\ref{X 2jump}) and recalling (\ref{delta}) we get the solution
\begin{eqnarray}\label{Xsol 2jump}
\hspace{-1cm}\left\{
\begin{array}{lll}
X_1(t)=c_1+(c_2+\mu_1)t+(c_3+\mu_2)\frac{t^2}{2}+\mu_3\frac{t^3}{6}+\sigma_1 W_1(t)\\
\quad \quad \quad \quad+\sigma_2\int_0^t(t-s)dW_2(s)+\sigma_3\int_0^t\frac{(t-s)^2}{2}dW_3(s)\\
\quad \quad \quad \quad +aH(t-\theta_1)+\frac{a}{\Delta}(t-\theta_0)I_{[\theta_0,\theta_1)}(t)\\
X_2(t)=c_2+(c_3+\mu_2)t+\mu_3\frac{t^2}{2}+\sigma_2 W_2(t)+\sigma_3\int_0^t(t-s)dW_3(s)\\
\quad \quad \quad \quad +\frac{a}{\Delta}I_{[\theta_0,\theta_1)}(t)\\
X_3(t)=c_3+\mu_3 t+\sigma_3 W_3(t)
\end{array}
\right.,
\end{eqnarray}
where $I_{[a,b)}(t)$ is the window function defined as
\begin{eqnarray}
I_{[a,b)}(t)=H(t-a)-H(t-b),
\end{eqnarray}
and the term $+aH(t -\theta_1)$ takes care of the phase deviation accumulated until time $\theta_1$.

The probability distribution function of the solution $\bf X$ at time $t$ is Normal with mean ${\bf M^{II}}(t)$ and covariance matrix ${\bf \Sigma}(t)$. Note that the mean can be decomposed in two terms where the first corresponds to the mean value of the process without anomalies (\ref{M}) 
\begin{eqnarray}\label{M II}
{\bf M^{II}}(t)={\bf M}(t)+\left[
\begin{array}{l}
aH(t-\theta_1)+\frac{a}{\Delta}(t-\theta_0)I_{[\theta_0,\theta_1)}(t)\\
\frac{a}{\Delta}I_{[\theta_0,\theta_1)}(t)\\
0
\end{array}
\right]
\end{eqnarray}
while its covariance matrix (\ref{Sigma}) does not change.

The solution in iterative form becomes:
\begin{eqnarray}\label{Xsol iter delta}
\hspace{-2cm}\left\{
\begin{array}{lll}
X_1(t_{k+1})=X_1(t_k)+(\mu_1+X_2(t_k))\tau+(\mu_2+X_3(t_k))\frac{\tau^2}{2}+\mu_3\frac{\tau^3}{6}+J_{k,1}\\
X_2(t_{k+1})=X_2(t_k)+(\mu_2+X_3(t_k))\tau+\mu_3\frac{\tau^2}{2}+J_{k,2}+\frac{a}{\Delta}\left[\delta_{(t_{k+1}-\theta_0)}-\delta_{(t_{k+1}-\theta_1)}\right]\\
X_3(t_{k+1})=X_3(t_k)+\mu_3 \tau+J_{k,3}
\end{array}
\right.,
\end{eqnarray}
where ${\bf J}_k=(J_{k,1},J_{k,2},J_{k,3})$ is (\ref{innovation}).
Again we note that the effect of the jump is present only in the affected component, $X_2$ in this case, and then the iterative process takes care of the effect on the other components.

\begin{remark}
Note that the algoritms can be easily extended to $n$ deterministic jumps.
\end{remark}

\subsection{Variance increase}
The atomic clock dynamics can also show changes in the variance. Let us model an atomic clock error with a change of the variance values in the interval $[\theta_0, \theta_1]$. 
The process for $t\notin [\theta_0, \theta_1]$ is described as the model without anomalies (\ref{X}) and its probability distribution is Normal with mean ${\bf M}(t)$ and covariance matrix ${\bf \Sigma}(t)$. At time $t\in [\theta_0, \theta_1]$ its distribution changes, while the mean ${\bf M}(t)$ is unchanged (\ref{M}) its covariance matrix becomes ${\bf \Sigma^{III}}(t)$ that is the same matrix (\ref{Sigma}) but with parameters  $\sigma_1$, $\sigma_2$ and $\sigma_3$ substituted by 
$\sigma'_1$, $\sigma'_2$ and $\sigma'_3$.

This anomaly can be modeled introducing a new innovation matrix in the time interval $[\theta_0, \theta_1]$.
The corresponding mathematical iterative model becomes
\begin{eqnarray}\label{Xsol iter var}
\hspace{-1cm}\left\{
\begin{array}{lll}
X_1(t_{k+1})=X_1(t_k)+(\mu_1+X_2(t_k))\tau+(\mu_2+X_3(t_k))\frac{\tau^2}{2}+\mu_3\frac{\tau^3}{6}+\tilde{J}_{k,1}\\
X_2(t_{k+1})=X_2(t_k)+(\mu_2+X_3(t_k))\tau+\mu_3\frac{\tau^2}{2}+\tilde{J}_{k,2}\\
X_3(t_{k+1})=X_3(t_k)+\mu_3 \tau+\tilde{J}_{k,3}
\end{array}
\right.,
\end{eqnarray}
where the innovation is

\begin{eqnarray}
{\bf \tilde{J}}_k=(\tilde{J}_{k,1},\tilde{J}_{k,2},\tilde{J}_{k,3})\sim 
\left\{
\begin{array}{ll}
N(0,{\bf Q'}) & \mbox{if } t_{k+1}\in[\theta_0,\theta_1]\\
N(0,{\bf Q}) & \mbox{if } t_{k+1}\notin[\theta_0,\theta_1]
\end{array}
\right.
\end{eqnarray}
where ${\bf Q}$ is the covariance matrix (\ref{Q}) while ${\bf Q'}$ is the same matrix where the parameters $\sigma_1$, $\sigma_2$ and $\sigma_3$ are substituted by 
$\sigma'_1$, $\sigma'_2$ and $\sigma'_3$.

\begin{remark}
In its present form the model accounts for anomalies at deterministic times. In reality anomalies occur at random times. It is possible to extend the model and the corresponding algorithms proposed to random times $\theta_i$, $i=0,1,2,3$ under the hypothesis that the random variables $\theta_i$ are independent for $i=0,1,2,3$ and that are independent from the process $\bf{X}$ without anomalies.
Under these hypothesis the distribution of the solution can be computed conditioning on the time of the anomalies. The iterative solution with deteriministic times can still be used for simulations after generating random times for the anomalies epochs according to their given distribution.
\end{remark}

\section{Simulations and examples}

The simulation of the models introduced in the previous paragraph can be easily handled using the iterative solutions (\ref{Xsol iter jump bis}), (\ref{Xsol iter delta}) and (\ref{Xsol iter var}).

In all the cases it is necessary to simulate the innovations ${\bf J}_k$ that, at any instant $t_k$, are independent identical distributed three dimensional Normal random variabiles with zero mean and covariance matrix  ${\bf Q}$ given by (\ref{Q}). Since the matrix ${\bf Q}$ is not diagonal, the three components of the vector ${\bf J}_k$ are dependent. Many softwares have routines to generate multidimensional Normal random variables. 

To better understand the correlation mechanism, we explicitely describe the procedure for generating ${\bf J}_k$ in our three dimensional case. 
It is well known that the vector ${\bf J}_k$ may be written as ${\bf J}_k={\bf AZ}$ where  ${\bf Z}$ is a standard Normal random vector with zero mean and identity covariance matrix ${\bf I}$, i.e. ${\bf Z}$ is a vector of three standard Normal independent random variables \cite{JP}. The matrix ${\bf A}$ transforms the three components of ${\bf Z}$ introducing a dependency that is captured in the correlation matrix ${\bf Q}$, where ${\bf Q=A A}^T$. Here the superscript $T$ denotes the transposed matrix.
The values of ${\bf A}$ can be computed using for example Cholesky decomposition method \cite{HJ} that gives 
\begin{eqnarray}
{\bf A}=\left[
\begin{array}{ccc}
\sqrt{q_{11}} & 0 & 0\\
 \frac{q_{21}}{\sqrt{q_{11}}} & \sqrt{\frac{q_{11}q_{22}-q_{21}^2}{q_{11}}}&0\\
 \frac{q_{31}}{\sqrt{q_{11}}} & \frac{q_{32}q_{11}-q_{31}q_{21}}{\sqrt{q_{11}(q_{11}q_{22}-q_{21}^2)}}  &\sqrt{q_{33}-\frac{q_{31}^2}{q_{11}}-\frac{q_{32}q_{11}-q_{31}q_{21}}{\sqrt{q_{11}(q_{11}q_{22}-q_{21}^2)}} } 
 \end{array}
\right]
\end{eqnarray}
where $q_{i,j}$, $i,j=1,2,3$ are the entry of ${\bf Q}$.

Note that the iterative solutions, useful for simulations, are exact solutions and do not introduce any approximation error. Exact simulations are helpful to understand the evolution of the clock behavior or to obtain estimates of functionals of clock states that cannot be computed theoretically.
The Matlab code can be freely downloaded from the INRIM website. 

 \subsection{Simulation of anomalous behavior}
Here we report some examples to illustrate the effect of the different types of anomalies of the process described in the previous section.

\begin{figure}[t!] 
\centering 
\includegraphics[scale=0.6]{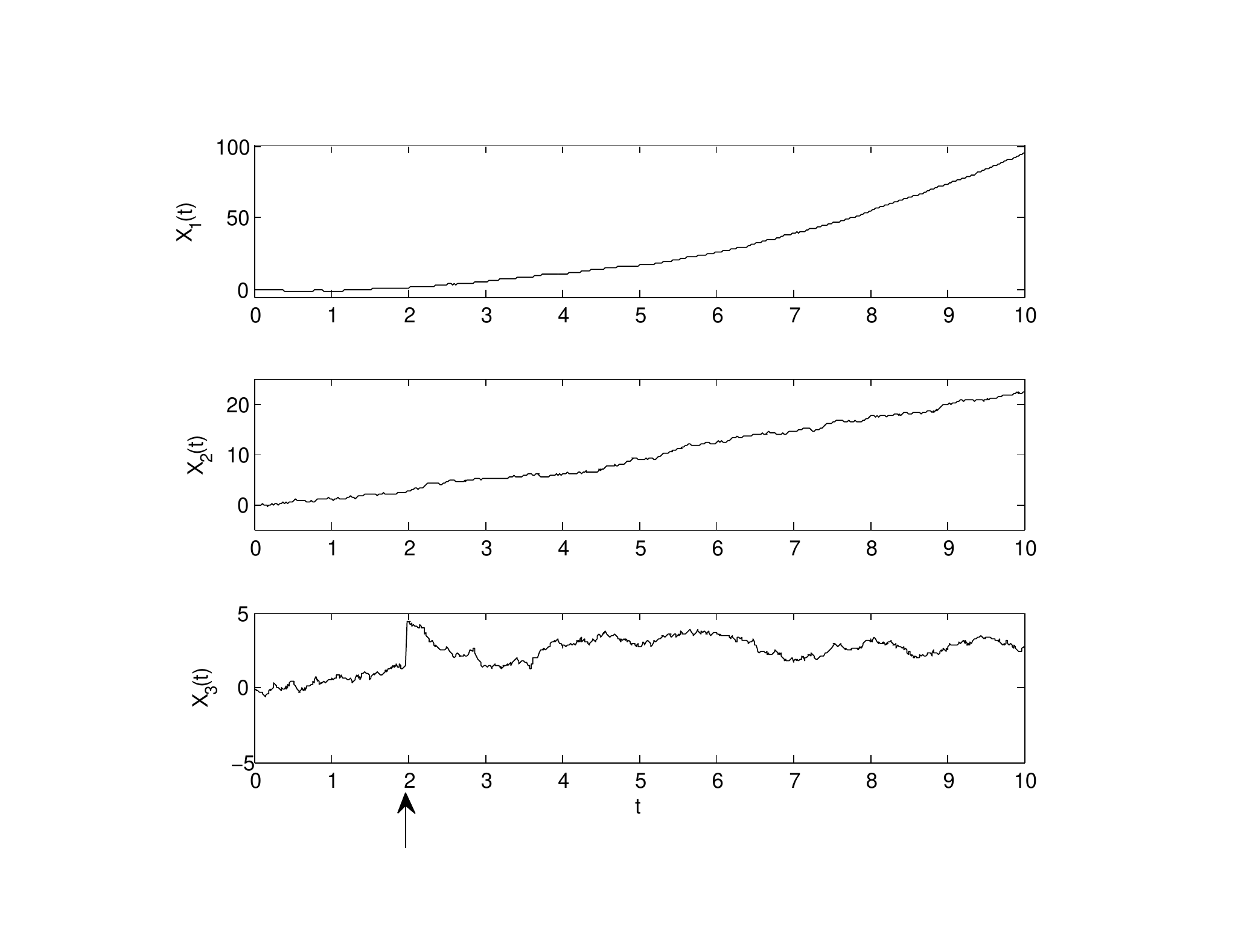}
\caption{Simulation of the three components of the process  $X(t)$ with parameters $\mu_1=\mu_2=\mu_3=0$,  $\sigma_1=\sigma_2=\sigma_3=1$, $a_1=a_2=0$, $a_3=3$ and $\theta_3=2$. The times of the jump is highlighted by an arrow (arbitrary units).}
\label{Fig:anom1 1 jump}
\end{figure}

\begin{figure}[t!] 
\centering 
\includegraphics[scale=0.6]{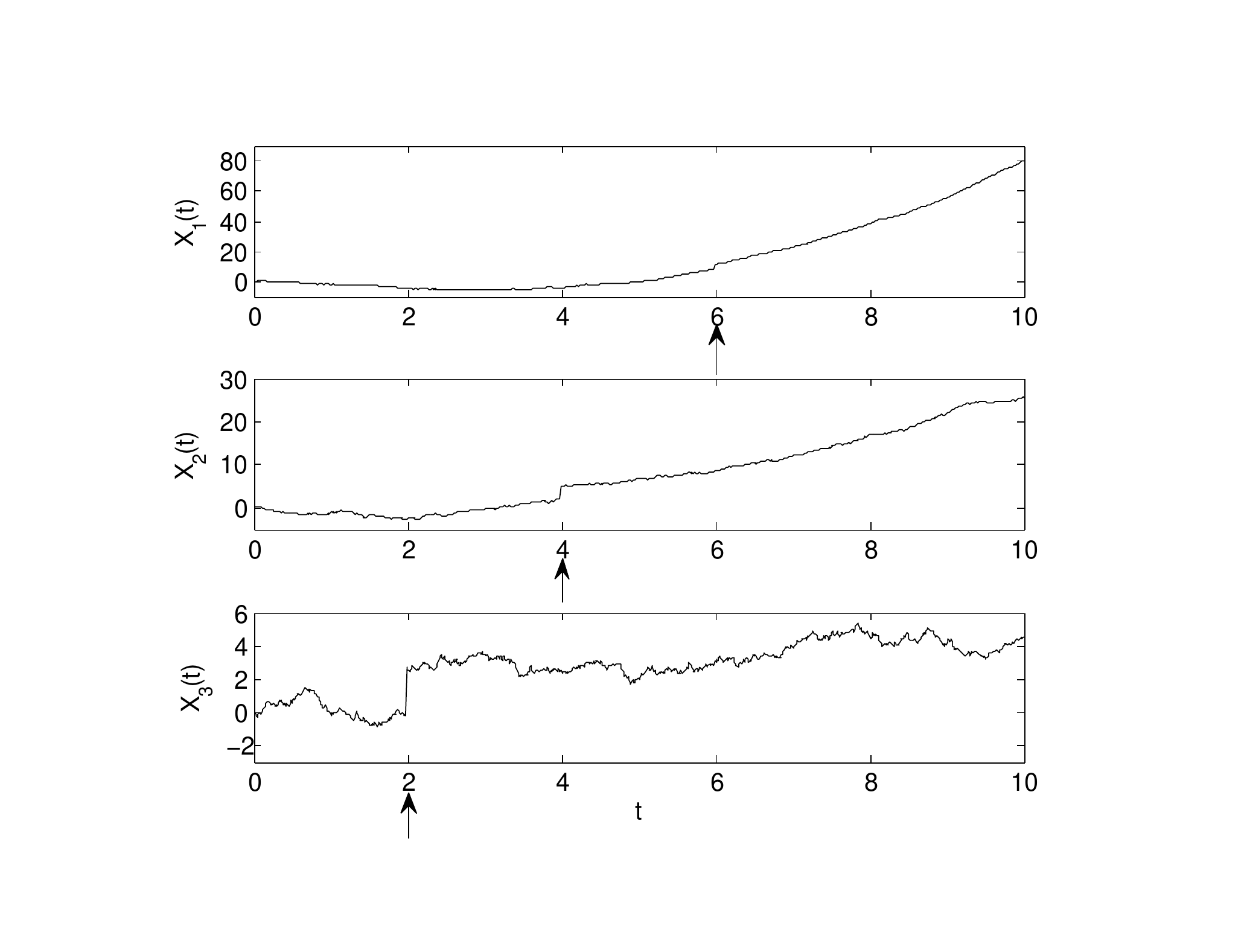}
\caption{Simulation of the three components of the process  $X(t)$ with parameters $\mu_1=\mu_2=\mu_3=0$,  $\sigma_1=\sigma_2=\sigma_3=1$, $a_1=a_2=a_3=3$ and $\theta_1=6$, $\theta_2=4$, $\theta_3=2$. The times of the jumps are highlighted by an arrow (arbitrary units).}
\label{Fig:anom1}
\end{figure}

In Figure \ref{Fig:anom1 1 jump} the model (\ref{X jump}) with an instantaneous jump on the third component is considered and a sample path of the process ${\bf X}(t)$ is plotted. The parameters considered are  $\mu_1=\mu_2=\mu_3=0$,  $\sigma_1=\sigma_2=\sigma_3=1$, $a_1=a_2=0$, $a_3=3$ and the times of the jump is $\theta_3=2$. The figure clearly shows the jump in the third component, the linear trend of the second component and the quadratic trend of the first one, as expected by the solution (\ref{Xsol jump}).

In Figure \ref{Fig:anom1} the model (\ref{X jump}) with instantaneous jumps in all the three components is considered and a sample path of the process ${\bf X}(t)$ is plotted. The parameters considered are  $\mu_1=\mu_2=\mu_3=0$,  $\sigma_1=\sigma_2=\sigma_3=1$, $a_1=a_2=a_3=3$ and the times of the jumps are $\theta_1=6$, $\theta_2=4$, $\theta_3=2$, respectively. The figure clearly shows the jumps in the trajectories and the increasing trend of each component sample path after the jumps. 

\begin{figure}[t!] 
\centering 
\includegraphics[scale=0.6]{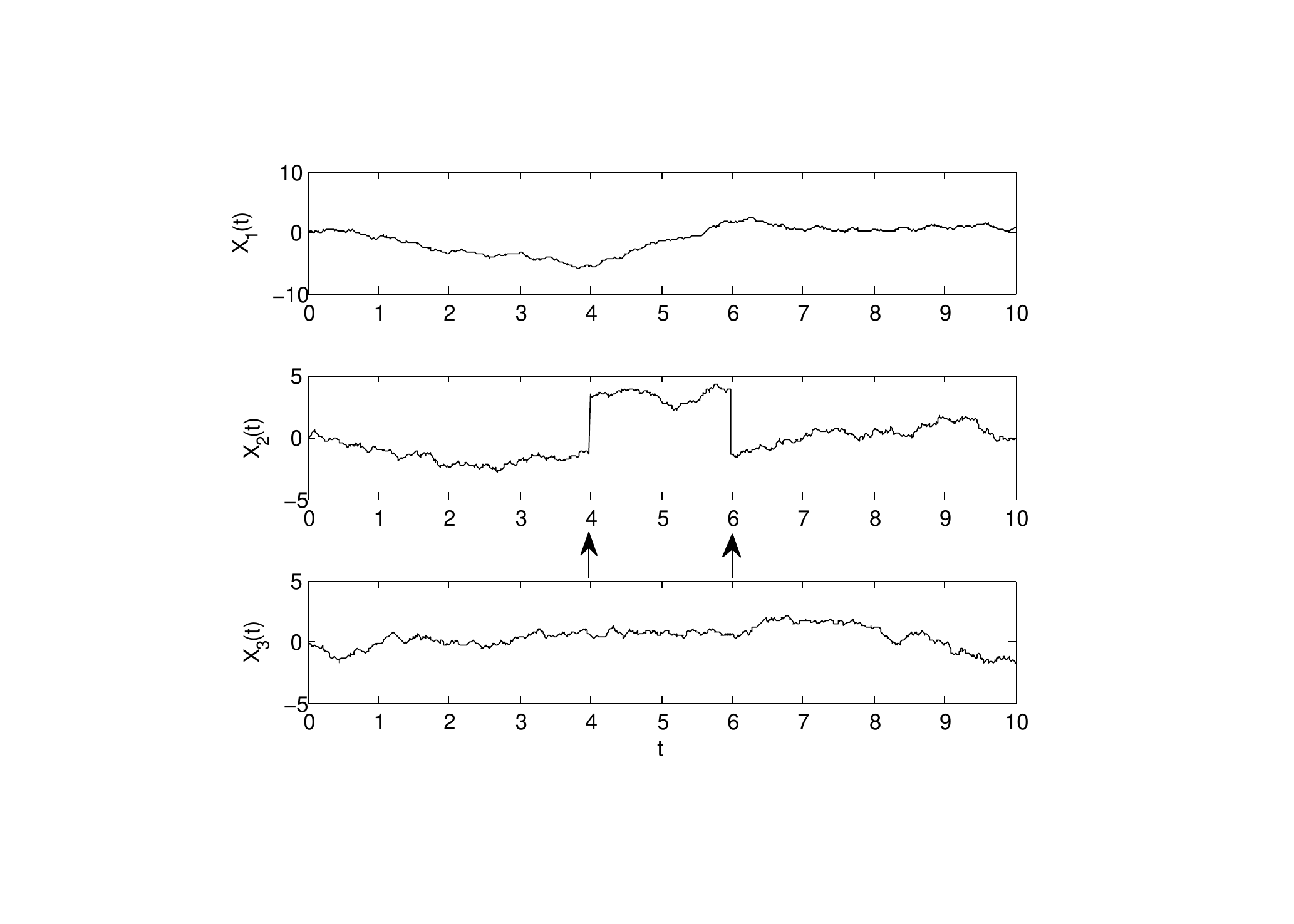}
\caption{Simulation of the three components of the process  $X(t)$ with parameters $\mu_1=\mu_2=\mu_3=0$,  $\sigma_1=\sigma_2=\sigma_3=1$, $a=4$ and the time interval of the jumps is $[4,6]$. The times of the jumps are highlighted by an arrow (arbitrary units).}
\label{Fig:anom2}
\end{figure}

In Figure \ref{Fig:anom2} the three components of the process $X(t)$ are plotted in case of two equal and opposite frequency jumps. The parameters considered are  $\mu_1=\mu_2=\mu_3=0$,  $\sigma_1=\sigma_2=\sigma_3=1$, $a=4$ and the time interval of the jumps is $[\theta_0, \theta_1]=[4,6]$. 
It is evident the linear trend of the first component in the interval $[4,6]$.

\begin{figure}[t!] 
\centering 
\includegraphics[scale=0.6]{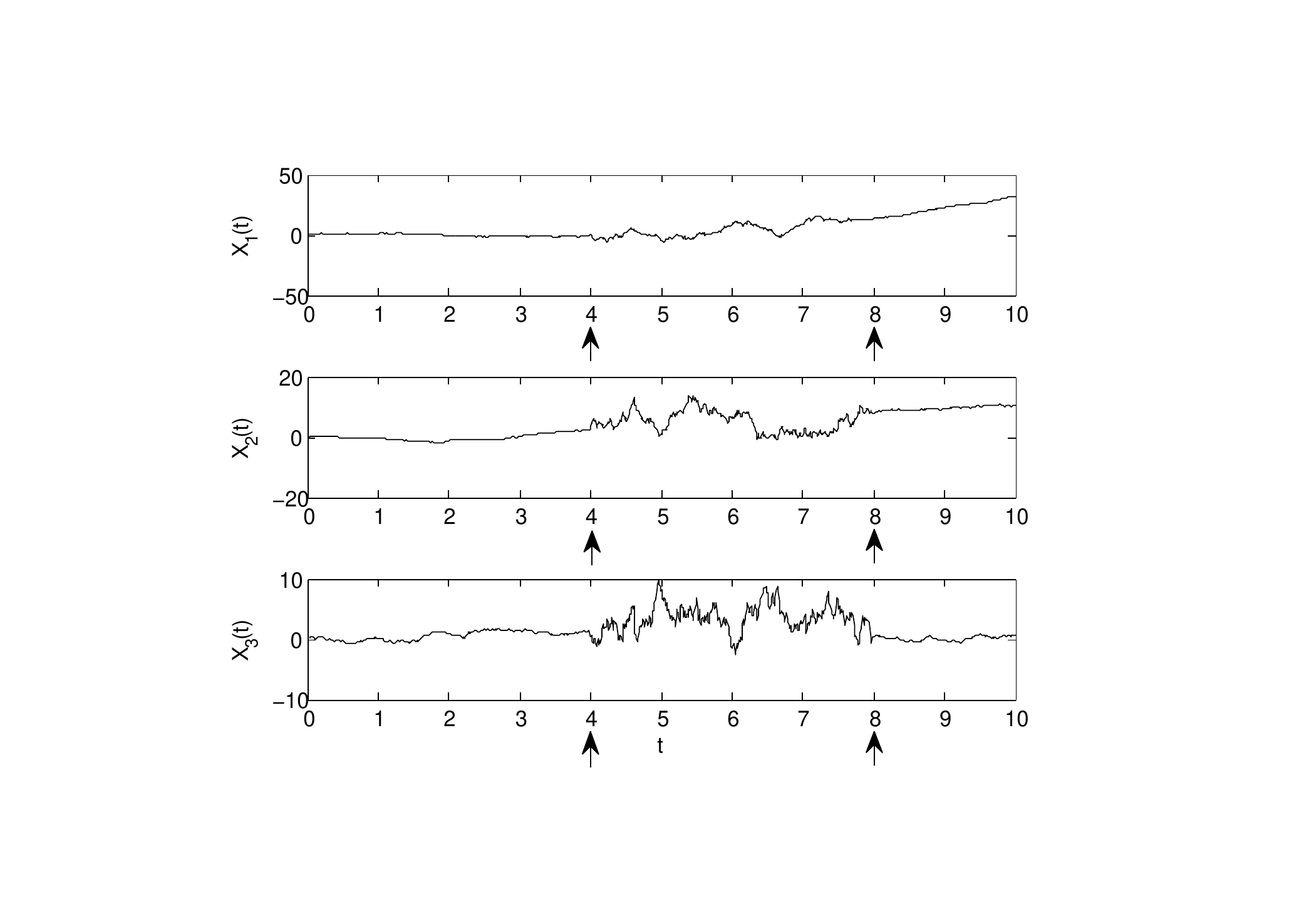}
\caption{Simulation of the three components of the process  $X(t)$ with parameters $\mu_1=\mu_2=\mu_3=0$,  $\sigma_1=\sigma_2=\sigma_3=1$,
$\sigma'_1=\sigma'_2=\sigma'_3=8$,  and the time interval of the change of variance is $[4,8]$. 
The time interval of variation of the variance is highlighted by arrows (arbitrary units).}
\label{Fig:anom3}
\end{figure}

In Figure \ref{Fig:anom3} the three components of the process $X(t)$ are plotted in case of variance increase. The parameters considered are  $\mu_1=\mu_2=\mu_3=0$,  $\sigma_1=\sigma_2=\sigma_3=1$,
$\sigma'_1=\sigma'_2=\sigma'_3=8$, 
and the time interval of the change of variance is $[4,8]$. The figure shows the effect of an increase of the noise in the selected interval.

The previous figures show that the model captures the main features of the evolution of the atomic clock error. 

\section{Application to a space GNSS clock}

In this section we tune the parameters of the model in order to study a case of interest concerning a Rubidium Atomic Frequency Standard (RAFS) placed on board of a satellite of a Global Navigation Satellite System (GNSS). We evaluate the effect of a clock anomaly on the clock prediction error, which is a crucial problem in navigation system.
As an example we consider the standard parameters describing the  behavior of a RAFS on board on the experimental satellites of the Galileo system \cite{WGBSHTT}.

The Allan deviation \cite{Al} indicating the stochastic short term noise of such RAFS is usually due to white frequency noise and it is assumed to be $\sigma_y(\tau)=5 \cdot 10^{-12}\tau^{-1/2}$.
The relationship between the Allan deviation and the diffusion coefficient of matrix $\bf Q$ is known \cite{ZT} and, since
\begin{equation}
\sigma_y(\tau) =\frac{\sigma_1}{\sqrt{\tau}},
\end{equation}
we can infer $\sigma_1=5\cdot 10^{-12}$. The other noises are negligible for this example, we set $\sigma_2= 10^{-22}$, $\sigma_3= 10^{-22}$ for numerical reasons only.  
We assume that all the polynomial, linear or quadratic slopes are zero, either because they are really negligible or because they have been estimated and appropriately compensated. We also assume that the clock is perfectly synchronized at epoch zero, therefore all the initial time constants are set to zero.

With this set of parameter we check which is the time deviation gained by such a clock a certain time after synchronization. Note that in the Galileo system the clocks are expected to be re-synchronized every 100 min = 6000 seconds \cite{WGBSHTT}.

Using the clock model without anomalies (\ref{X}) the time deviation $X_1$ at time $t= 6000$ s is a Normal random variable with mean value ${\bf M}_1(t)= E[X_1(t)]= 0$ and standard deviation approximately $\sqrt{\bf\Sigma_{11}(t)}\simeq \sigma_1 \sqrt{t}= 5\cdot10^{-12} \sqrt{6000}\simeq 4\cdot 10^{-10}$ s $= 0.4$ ns.
In this case a time deviation interval of $\pm 1,96 \sqrt{{\bf \Sigma}_{11}(t)}\simeq \pm759\cdot 10^{-12}$ s $=\pm 0.8$ ns, with 95\% confidence level, after 6000 s, represents the prediction error and hence the contribution on the positioning error due to unpredictable clock time deviation. Note that the polynomial behavior of the three components has not been considered in this example.

\begin{figure}[t!] 
\centering 
\includegraphics[scale=0.6]{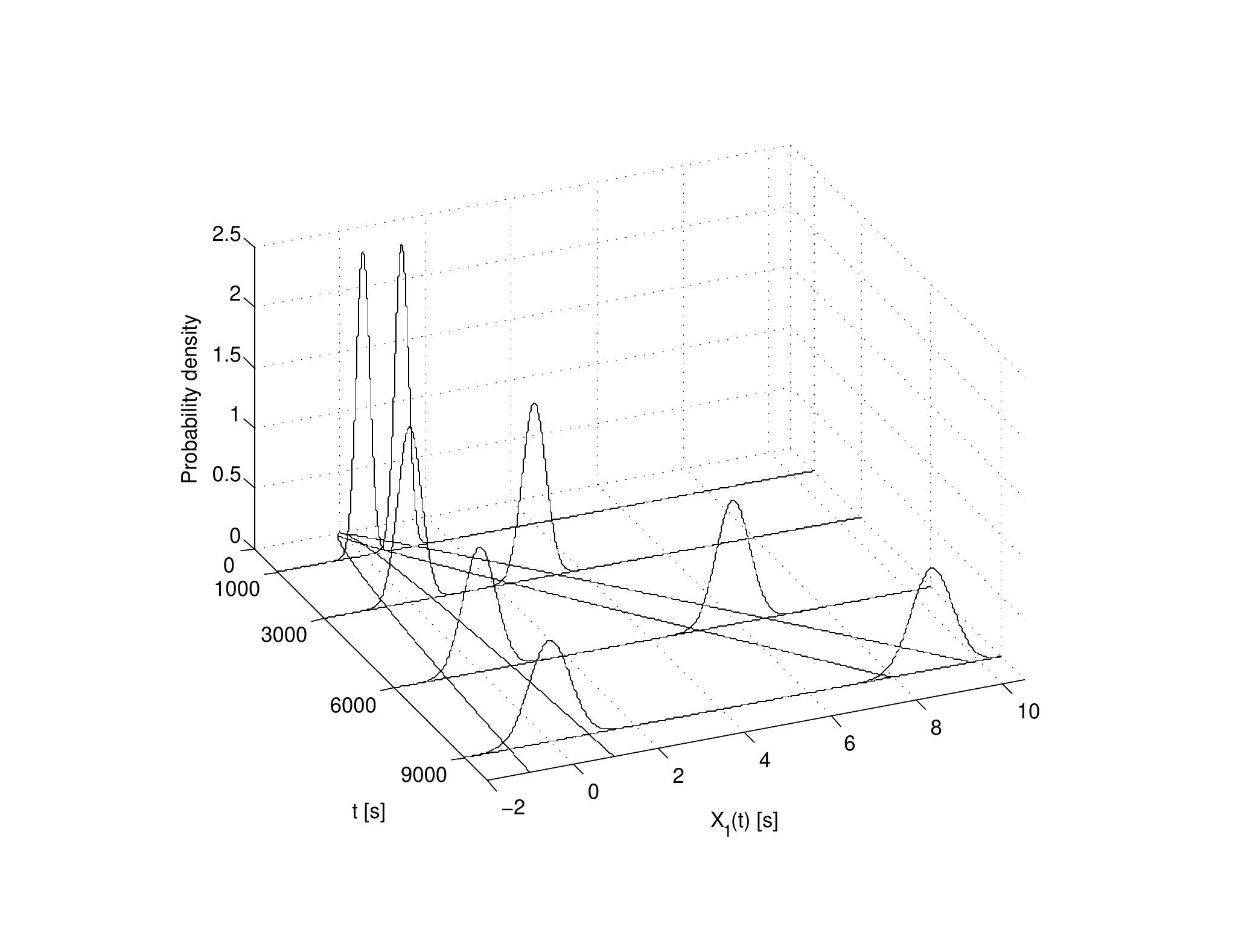}
\caption{Probability density function and the 95\% confidence interval of the time deviation $X_1$ at time $t=[1000, 3000,6000,9000]$ s, with parameters $\mu_1=\mu_2=\mu_3=0$,  $\sigma_1=5\cdot 10^{-12}$, $\sigma_2=\sigma_3=10^{-22}$ and jump $a_2=10^{-12}$ at time $\theta=100$ compared with the case without anomalies.}
\label{Fig:IdC_gauss_3D_bis}
\end{figure}

To check the effect of a clock anomaly, we consider a frequency jump of relative value $10^{-12}$. 
Frequency jumps are sometimes observed on space clock on board of navigation satellites \cite{CFSGT}. In our notation a frequency jump corresponds to (\ref{Xsol jump}) when we choose $a_2=10^{-12}$ and  $a_1=a_3=0$. We suppose that the frequency jump occurs at epoch $\theta$ and the effect on the time deviation at epoch $t$ depends on the interval $(t-\theta)$. Looking to (\ref{Xsol jump}) we expect the time deviation $X_1$ to be a Normal random variable with mean value ${\bf M^I}_1(t)=a_2 (t-\theta)$ and standard deviation $\sqrt{\bf\Sigma_{11}(t)}= \sigma_1 \sqrt{t}= 5\cdot10^{-12} \sqrt{t}$.
Let us deal with anomaly occurring at $\theta = 100$ s or at $\theta = 5000$ s. In the first case the anomaly occurs immediately after the synchronization and at $t= 6000$ s its effect is quite impacting. The 95\% confidence level for the time deviation $X_1$ gives the interval $a_2 (t-\theta)\pm 1.96 \sqrt{\Sigma_{11}(t)}= 5.9 \pm 0.8$ ns. 

In the second case, the anomaly occurs just at the end of the predictive period, and the prediction error at 95\% confidence level in that case is $a_2 (t-\theta)\pm 1.96 \sqrt{\Sigma_{11}(t)}=1 \pm 0.8$ ns, 
therefore the effect is much smaller, very close to the case without anomalies discussed above. 

In Figure \ref{Fig:IdC_gauss_3D_bis} is shown the probability density function and the 95\% confidence interval of the time deviation $X_1$ at time $t=[1000, 3000,6000,9000]$ s, with jump $a_2=10^{-12}$ at time $\theta=100$ compared with the case without anomalies.

As we see from (\ref{M I}), the impact of the frequency jump on the time deviation and hence on the time prediction error is linearly depending on $t-\theta$, the worst case being at $\theta=0$ giving a confidence interval $a_2 t \pm 1.96 \sqrt{\Sigma_{11}(t)}$, while the best case of no anomaly gives $\pm 1.96 \sqrt{\Sigma_{11}(t)}$
at 95\% confidence level.

\section{Conclusion}
A previous mathematical model based on multidimensional stochastic differential equation describing the atomic clock error has been extended to include the presence of clock anomalies experimentally observed on clocks, particularly in GNSS space clocks. The model includes phase, frequency, and drift jumps, as well as changes on the random noise power. Examples of such trends are presented and the case of a GNSS space clock is discussed by evaluating the time deviation gained at a certain time after synchronization with or without the effect of a frequency jump. Our results show the ability of the proposed model to describe clock anomalies and their effect on complex systems.

\section*{Acknowledgment}
Work partially supported by Turin University: Stochastic Processes/ ZUCC01CT11.

\end{document}